\documentclass[12pt]{article}
\usepackage{amsmath,amsfonts,amssymb,amsthm}

\newcommand{\al}{\alpha}
\newcommand{\del}{\delta}					

\newcommand{\W}{\Omega}
\newcommand{\lm}{\lambda} 
\renewcommand{\th}{\theta}
\newcommand{\w}{\omega}
\newcommand{\1}{{\bf 1}}
\newcommand{\id}{{\rm id}}

\newcommand{\End}{{\rm End\,}}
\newcommand{\wt}{{\rm wt}}
\newcommand{\Res}{{\rm Res\,}}

\def\Z{\mathbb{Z}}

\def\C{\mathbb{C}}
\def\N{\mathbb{N}}
\def\h{\mathfrak{h}}
\def\g{\mathfrak{g}}

\def\hh{\hat{\h}}

\newcommand{\NO}{\,{\raise0.25em\hbox
{$\mathop{\hphantom {\cdot}}
\limits^{_{\circ}}_{^{\circ}}$}}\,}
\newcommand{\M}[1]{M(1)^{#1}}
\newcommand{\Ml}[1]{M(1,#1)}
\newcommand{\V}[1]{V_{L}^{#1}}
\newcommand{\Vl}[1]{V_{#1 +L}}

\newcommand{\Va}[1]{V_{\al/2+L}^{#1}}

\newcommand{\ots}{\otimes}
\newcommand{\ops}{\oplus}

\newcommand{\eqa}{\begin{eqnarray}}
\newcommand{\eeqa}{\end{eqnarray}}
\newcommand{\eqn}{\begin{eqnarray*}}
\newcommand{\eeqn}{\end{eqnarray*}}
\newcommand{\nn}{\nonumber}

\newtheorem{dfn}{Definition}[section]
\newtheorem{pro}[dfn]{Proposition}
\newtheorem{thm}[dfn]{Theorem}

\newtheorem{lem}[dfn]{Lemma}
\newtheorem{cor}[dfn]{Corollary}
\newtheorem{rem}[dfn]{Remark}

\makeatletter
\@addtoreset{equation}{section}

\makeatother
\def\bl{\begin{lem}\label}
\def\el{\end{lem}}
\def\bt{\begin{thm}\label}
\def\et{\end{thm}}
\def\bp{\begin{pro}\label}
\def\ep{\end{pro}}
\def\br{\begin{rem}\label}
\def\er{\end{rem}}
\def\bc{\begin{cor}\label}
\def\ec{\end{cor}}
\def\bd{\begin{dfn}\label}
\def\ed{\end{dfn}}
\def\proof{{\it Proof. }}
\def\qed{\hspace*{\fill}{\hbox{ $\square$}}\vspace{5mm}}
\newcommand{\Pv}[1]{\mathcal{P}({#1})}
\def\Av{A(V)}
\newcommand{\Com}[2]{\mbox{$\left(\begin{array}{c}#1\\#2
\end{array}\right)$}}
\newcommand{\T}[1]{\tilde{#1}}
\def\B{\langle}
\def\K{\rangle}
\def\span{{\rm span}}

\begin{document}
\title{\Large  The charge conjugation orbifold $V_{\Z\al}^{+}$ is
rational\\ when ${\B\al,\al\K}/{2}$ is prime}
\author{\large Toshiyuki Abe}
\date{\small\it Department of Mathematics, Graduate School of
Science, Osaka University\\
 Toyonaka, Osaka 560-0043, Japan\\
{\rm e-mail: sm3002at@ecs.cmc.osaka-u.ac.jp}}
\maketitle
 
\begin{abstract}
In this paper we prove the rationality of the vertex operator
algebra
$\V{+}$ for a rank one positive definite even lattice $L$ whose
generator has a prime half length. 
\end{abstract}

\baselineskip 5mm


\section{Introduction}
Let $L=\Z\al$ be a rank one even lattice defined by $\B \al,\al\K=2k$ for
a positive integer $k$. Then it is known that  the lattice vertex operator
algebra $\V{}$ is rational (see \cite{D1}). The vertex operator algebra 
$\V{}$ has an automorphism induced from the
$-1$-isometry of $L$, and the fixed point set $\V{+}$ is a simple vertex
operator algebra (cf. \cite{DN2}), which is called the charge conjugation
orbifold in the literature \cite{DVVV} and 
\cite{KT}. In this paper we show that the vertex operator algebra
$\V{+}$ is rational when $k$ is a {\it prime integer}. 

The rationality of $\V{+}$ with $k=3$ may help to study the
construction of moonshine module. The moonshine vertex operator
algebra $V^{\natural}$ has a vertex operator subalgebra isomorphic to the
tensor product
$U=(\V{+})^{\ots24}$ of $\V{+}$ with $k=3$ which has the same Virasoro
element of
$V^{\natural}$ (\cite{Sh}).  As
proved in \cite{FHL}, the tensor product $V=\ots_{i=1}^{n}V^{i}$ of rational
vertex operator algebras $V^{i}$ is also rational and any irreducible
$V$-modules are isomorphic to tensor products $\ots_{i=1}^{n}M^{i}$ of
some irreducible $V^{i}$-module $W^{i}$.  Therefore by the rationality of
$\V{+}$ with $k=3$, one knows that $V^{\natural}$ is decomposed into a
direct sum of irreducible $U$-modules which are tensor products of 24
irreducible $\V{+}$-modules. The similar decompositions of
$V^{\natural}$ as a direct sum of irreducible modules of the tensor
product $L(1/2,0)^{\ots48}$ of the Virasoro vertex operator algebra
$L(1/2,0)$ are known (cf. \cite{DMZ} and \cite{DGH}), and these
decompositions are used to compute the characters of elements of the
Monster simple group (see \cite{M5}). More precisely, Miyamoto
\cite{M2} constructed a vertex operator algebra $M_{D}$ associated to an
even linear binary code $D$ in $\Z_{2}^{n}$ for
$n\in\N$, which is called a code vertex operator algebra. The vertex 
operator algebra $M_{D}$ contains $L(1/2,0)^{\ots n}$ as a vertex
operator subalgebra, and each irreducible components of $M_{D}$ as
an $L(1/2,0)^{\ots n}$-module corresponds to a codeword of
$D$. In
\cite{M5}, the moonshine vertex operator algebra is reconstructed as an
induced module, which is introduced and studied in \cite{DLi}, of some
code vertex operator algebra, and the characters of $2A$, $2B$ and $3C$
elements of the Monster are calculated explicitly. On the other hand, in
\cite{Sh} the Mckay-Thompson series for
$4A$ elements of the Monster are calculated by using a decomposition of
$V^{\natural}$ as a
$(\V{+})^{\ots 24}$-module. Thus this decomposition 
may shed light on another aspects of the actions of the Monster on the
moonshine vertex operator algebra.  

In \cite{M3}, the structures of modules for code vertex
operator algebras are studied in detail, and there the rationality and
fusion rules of $L(1/2,0)$ play an important role. The fusion rules for
charge conjugation orbifolds
$\V{+}$ are completely determined in
\cite{A}. Together with the rationality of $\V{+}$ with
$k=3$, we hope that the analogue of the representation theory for code
vertex operator algebra can be developed for the vertex operator algebra 
$\V{+}$ with
$k=3$.  For example, as been studying by Lam and Shimakura (\cite{LS}),
the induced module construction from
$\V{+}$ may be explained by means of the language of code. 

Let $V$ be a vertex operator algebra. In \cite{Z}, an
associative algebra $A(V)$ associated with $V$ is introduced.  Zhu's
algebra $A(V)$ has many information on modules for $V$, for example,
there is a bijection between the set of equivalence classes of irreducible 
$V$-modules and that of irreducible $A(V)$-modules. Furthermore if 
$V$ is rational, then Zhu's algebra $A(V)$ is  semisimple.

Suppose that Zhu's algebra $A(V)$ is semisimple. Then we can prove
that any admissible
$V$-module $M$ has a generalized eigenspace decomposition for
$L(0)=\w(1)$, where $\w$ is the Virasoro element of $V$ (see
Theorem \ref{MP1}). This follows from the fact that if $A(V)$ is
semisimple, then the lowest degree space of an admissible $V$-module
$M$ is a direct sum of irreducible
$A(V)$-modules and then $L(0)$ acts semisimply on  the $V$-submodule
generated by the lowest degree space. Let
$\mathcal{P}(V)$ be the set of lowest weights of all irreducible
admissible 
$V$-modules. Then we prove that for any admissible $V$-module
and any complex number $\lm$, the direct sum
$\ops_{n\in\Z}M^{(\lm+n)}$ of generalized eigenspace $M^{(\lm+n)}$ of
eigenvalue $\lm+n$ $(n\in\Z)$ is an admissible $V$-submodule
of $M$. We also show that if $\mu\not\in\Pv{V}+\N$ then the
generalized eigenspaces $M^{(\mu)}$ is zero. These imply that $M$ is a
direct sum of admissible submodules of the form
$\ops_{n=0}^{\infty}M^{(\lm+n)}$ for some $\lm\in\mathcal{P}(V)$.  Thus
in order to prove that $V$ is rational, it suffices to show that for any
$\lm\in\Pv{V}$ every admissible $V$-module of the form
$\ops_{n=0}^{\infty}M^{(\lm+n)}$ are completely reducible. If
$\lm\in\mathcal{P}(V)$ satisfies the 
condition $(\lm+\Z_{+})\cap\Pv{V}=\emptyset$, then it is not hard to
see that any admissible modules $M=\ops_{n=0}^{\infty}M^{(\lm+n)}$ are 
completely reducible. Hence if we show that for any
$\lm\in\mathcal{P}(V)$ not
satisfying $(\lm+\Z_{+})\cap\Pv{V}=\emptyset$, every admissible
$V$-module of the form $M=\ops_{n=0}^{\infty}M^{(\lm+n)}$ is
completely reducible, then we have the rationality of $V$. 
  
Since Zhu's algebra $A(\V{+})$ is semisimple (see \cite{DN2}), we
can use the results explained in the previous paragraph to prove the
rationality of the vertex operator algebra $\V{+}$. The classification of
irreducible $\V{+}$-modules in \cite{DN2} shows that the set of lowest
weights for
$\V{+}$ is $\mathcal{P}(\V{+})=\{\,0,1,r^{2}/4k, 1/16,9/16\,\}$. Under the
assumption that $k$ is prime, all numbers
$0,1,1/16,9/16$ and $r^{2}/4k$ for $1\leq r\leq k$ are distinct and that 
lowest weights
$\lm$ except $0$ satisfy the condition 
$(\lm+\Z_{+})\cap\Pv{\V{+}}=\emptyset$ (see Lemma \ref{Lem4}).
Thus to show that $\V{+}$ with a prime integer $k$ is rational, it
suffices to prove that every admissible $\V{+}$-module of the form
$M=\ops_{n\in\N}M^{(n)}$ is completely reducible. We will prove
that the submodules of $M$ generated by $M^{(0)}$ and $M^{(1)}$ are
completely reducible because we can easily see that $M$ is generated by
$M^{(0)}$ and $M^{(1)}$.  

The organization of this paper is as follows: In Section \ref{S2.1}, we
review definitions of modules for a vertex operator algebra, and recall
some results related to Zhu's algebra in Section \ref{S2.2}.
In Section \ref{S3}, we prove that for a vertex operator algebra
$V$ whose Zhu's algebra is semisimple, every admissible module is a
direct sum of generalized eigenspaces for $L(0)$. The constructions
of the vertex operator algebra $\V{+}$ and its irreducible modules are
explained in Section \ref{S4.1}, and in Section \ref{S4.2} we prove that
the rationality of $\V{+}$ with a prime integer $k=\B\al,\al\K/2$ by using
results in Section \ref{S3}.  

Throughout the paper, $\N$ is the
set of nonnegative integers and $\Z_{+}$ is the set of positive integers.

\section{Preliminaries}\label{S2}
In Section \ref{S2.1}, we recall the definition of 
modules over a vertex operator algebra,  following
\cite{FLM}, \cite{FHL}. We also review the definition of  Zhu's algebra and
some related results in Section \ref{S2.2}.


\subsection{Modules for a Vertex Operator Algebra}\label{S2.1}

Let $(V,Y,\1,\w)$ be a vertex operator algebra (cf. \cite{FLM} and
\cite{FHL}). We write $Y(a,z)=\sum_{n\in\Z}a(n)z^{-n-1}$ for $a\in V$.
Then $V$ has an eigenspace decomposition $V=\ops_{n\in\Z}V(n)$ for
$L(0)=\w(1)$, where $V(n)$ is the eigenspace of eigenvalue $n$ for
$n\in\Z$. A vector $a$ in $V(n)$ is called to be homogeneous of weight
$n$. We denote the weight of $a\in V(n)$ by $n=\wt(a)$. 

\bd{Def1} A {\rm weak $V$-module} is a pair $(M,Y_{M})$ of
a vector space $M$ and a linear map
\eqn
Y_{M}:V&\to&(\End M)[[z,z^{-1}]],\\
a&\mapsto&Y_{M}(a,z)=\sum_{n\in\Z} a(n) z^{-n-1},\
(a(n)\in \End M)
\eeqn
such that the following conditions are satisfied for $a,b\in V$ and
$u\in M$:

\noindent
$(1)$
\eqa
Y_{M}(a,z)v\in M((z)),\label{GC}
\eeqa
\noindent
$(2)$ {\rm (The Jacobi identity).}
\eqa
&&z_{0}^{-1}\delta\left(
{\frac{z_{1}-z_{2}}{z_{0}}}\right)
Y_{M}(a,z_{1})Y_{M}(b,z_{2})-z_{0}^{-1}\delta
\left({\frac{z_{2}-z_{1}}{
-z_{0}}}\right)Y_{M}(b,z_{2})Y_{M}(a,z_{1})\nn\\
&&{ }=z_{2}^{-1}\delta\left( {\frac{z_{1}-z_{0}}{z_{2}}}\right)Y_{M}
(Y(a,z_{0})b,z_{2}),\label{Jac}
\eeqa

\noindent
$(3)$ 
$$ Y_{M}(\1,z)=\id_{M}.$$
\ed

Set $Y_{M}(\w,z)=\sum_{n\in\Z}L(n)z^{-n-2}$. Then as proved in
\cite{DLM1}, one have the commutation relations of the Virasoro algebra
$$[L(m),L(n)]=(m-n)L(m+n)+\frac{m^{3}-m}{12}\delta_{m,-n}\,
c_{V}\id_{M},$$
for $m,n\in\Z$, where $c_{V}$ is the central charge of $V$.  We also have
the
$L(-1)$-derivative property
\eqa 
Y_{M}(L(-1)a,z)=\frac{d}{dz}Y(a,z)\hbox{, for all $a\in V$}.\label{DP1}
\eeqa
A weak $V$-module $(M,Y_{M})$ is often denoted by $M$. 
A {\it weak $V$-submodule} of a weak module
$M$ is a subspace $N$ of $M$ such that $a(n)N\subset N$ hold for all 
$a\in V$ and $n\in\Z$. 

By the Jacobi identity (\ref{Jac}), we
have the following commutation relations for
$a,b\in V$ and $m,n\in\Z$:
\eqa
[a(m),b(n)]=\sum_{i=0}^{\infty}\Com{m}{k}(a(k)b)(m+n-k)\label{CR}
\eeqa 

Now we set $\T{a}(n)=a(\wt(a)+n-1)$ for a homogeneous vector 
$a\in V$ and $n\in \Z$, and extend it to every $a\in V$ by linearity.
Set 
$$\W(M)=\{\, u\in M\,|\, \T{a}(n)u=0\hbox{ for any $a\in V$ and
$n\in\Z_{+}$}\,\}.$$ 
By the commutation relation (\ref{CR}), we see that
$\W(M)$ is preserved  under the action of $o(a)=\T{a}(0)$ for $a\in V$.
  
For a weak $V$-module $M$ and $u\in M$, we denote by $\B u\K$ the
weak submodule of $M$ generated by $u$. Then
the following lemma holds: 


\bl{Lem1} {\rm (\cite{L1})} Let $M$ be a weak $V$-module and $u\in M$.
Then 
$$\B u\K=\span_{\C}\{ a(n)u\,|\,a\in V\hbox{ and }n\in\N\,\}.$$
\el


\bd{Def2}An {\rm admissible $V$-module} $M$ is a weak 
$V$-module equipped with an $\N$-grading $M=\ops_{n=0}^{\infty}
M_{n}$ such that the condition 
\eqa
a(n)M_{m}\subset M_{\wt(a)+m-n-1}\label{AD1}
\eeqa
holds for any homogeneous vector $a\in V$ and $m,n\in\Z$, where
$M_{m}=0$ for $n<0$.
\ed
For an admissible $V$-module $M$, we
have the linear map $d:M\to M$ defined by
$d|_{M_{n}}=n\,\id_{M_{n}}$ ($n\in\N$), which is called the degree
operator of $M$.  By (\ref{AD1}), we have commutation relation
\eqa [d,\T{a}(n)]=-n\,\T{a}(n) \label{CR2}
\eeqa
for $a\in V$ and $n\in\Z$. 

A $d$-invariant weak $V$-submodule $N$ of an admissible $V$-module
$M$ is called an {\it admissible
$V$-submodule} of $M$. Then $N$ is decomposed into a direct
sum of eigenspaces for $d$ as
$N=\ops_{n\in\N}N_{n}$, where $N_{n}=N\cap M_{n}$ for $n\in \N$. For
an admissible $V$-module $M$ and its admissible $V$-submodule $N$,
the quotient space
$M/N$ has a natural admissible $V$-module structure; the homogeneous
space of degree $n$ is given by 
\eqa
(M/N)_{n}=(M_{n}+N)/N.
\label{QM}
\eeqa

Now we recall the definition of the rationality of a vertex operator
algebra.  An {\it irreducible} admissible $V$-module 
$M$ is an admissible $V$-module which has no  admissible
$V$-submodule except $0$ and $M$. When $M$ is a direct sum of
irreducible admissible submodules, $M$ is called to be completely
reducible.
\bd{Def4}
A vertex operator
algebra
$V$ is said to be {\it rational} if any admissible $V$-module is completely
reducible. 
\ed


\subsection{ Zhu's Algebra $\Av$}\label{S2.2}
In this subsection we review the definition of Zhu's algebra $\Av$
associated to a vertex operator algebra $V$ and recall some related
results in \cite{Z} and \cite{DLM1}. 

Let $(V,Y,\1,\w)$ be a vertex operator algebra. For homogeneous vectors
$a\in V$ and $b\in V$, we define 
\eqa 
a\circ b=\Res_{z}\frac{(1+z)^{\wt(a)}}{z^{2}}Y(a,z)b.
\eeqa
Let $O(V)$ be the subspace of $V$ spanned by all vectors $a\circ b$
for homogeneous vectors $a,b\in V$, and set $\Av=V/O(V)$. We denote
the image of $a\in V$ in $A(V)$ by $[a]$.  For any homogeneous vector
$a\in V$ and $b\in V$, define  
\eqn 
a*b=\Res_{z}\frac{(1+z)^{\wt(a)}}{z}Y(a,z)b
\eeqn
and extend it on V for $a$ by linearity.  The following
theorem is due to \cite{Z} and \cite{DLM1}:


\bt{T2} (1) The operation $*$
induces an associative algebra structure on $\Av$, and the image
$[\w]$ is contained in the center of $\Av$.

(2) For any admissible $V$-module $M$, the linear
map $V\to\End \W(M)$ defined by $a\mapsto o(a)|_{\W(M)}$ for $a\in
V$ induces an $\Av$-module structure on $\W(M)$. 

(3) The map $M\mapsto\W(M)$ induces a bijection from
the set of all equivalence classes of irreducible admissible  
$V$-modules to the set of all equivalence classes of irreducible
$\Av$-modules.   

(4) If $A(V)$ is semisimple, then there are only finitely many
inequivalent irreducible admissible $V$-modules. 
\et

It is known that for any irreducible admissible $V$-module
$M$, there exists a complex number $\lm\in\C$ such that $M$ is
decomposed into a direct sum of $L(0)$-eigenspaces as 
$M=\ops_{n=0}^{\infty}M(\lm+n)$, where
$M(\lm+n)$ is the $L(0)$-eigenspace of eigenvalue $\lm+n$ (see
\cite[Lemma 1.2.1]{Z}). We call
the complex number $\lm$ a {\it lowest weight} of $M$. We denote by
$\mathcal{P}(V)$ the set of lowest weights of all irreducible admissible
$V$-modules. 

\bp{Pro1} (1) For an admissible $V$-module $M$, if a
subspace $W$ of $\W(M)$ is invariant under the action of $o(a)$ for any
$a\in V$, then $W$ is $\Av$-submodule of $\W(M)$.

(2) Let $W$ be a completely reducible $\Av$-module. Then
every eigenvalues for $[\w]$ on $W$ are in $\mathcal{P}(V)$. 
\ep
\proof   The assertion (1) follows from Theorem \ref{T2} (2). We
shall prove (2). Since W is completely reducible as an $\Av$-module, it is
decomposed into a direct sum of irreducible $\Av$-submodules. Since
$[\w]$ is in the center of $\Av$, $[\w]$ acts on each irreducible
component by a scalar. By Theorem \ref{T2} (3), any irreducible
$\Av$-module is isomorphic to $\W(M)$ for some irreducible
admissible $V$-module $M$. Since $[\w]$ acts on $\W(M)$ as $L(0)$, we
see that the eigenvalues for
$[w]$ on irreducible components of $W$ are
lowest weights. This proves (2). \qed


\section{Generalized Eigenspaces Decompositions\\ of Admissible
$V$-Modules}\label{S3}
Let $V$ be a vertex operator algebra. For a weak $V$-module $M$,
we denote the generalized eigenspace for $L(0)$ with eigenvalue
$\lm\in\C$ by
$M^{(\lm)}$, i.e., 
$$M^{(\lm)}=\{\,u\in M\,|\,(L(0)-\lm)^{k}u=0\hbox{ for some
$k\in\N$}\,\}.$$

The purpose of this section is to show the following theorem:
\bt{MT} Let $V$ be a vertex operator algebra. Suppose that Zhu's
algebra $\Av$ is semisimple. Then any admissible
$V$-module
$M$ is decomposed into a direct sum of generalized eigenspaces for
$L(0)$ as follows:
$$M=\ops_{i=1}^{s}(\ops_{n\in\N}M^{(\lm_{i}+n)})$$
for some $\lm_{i}\in\mathcal{P}(V)$ $(i=1,2,\cdots, s)$. Moreover, for
each $\lm_{i}$, the subspace $\ops_{n\in\N}M^{(\lm_{i}+n)}$ is an
admissible $V$-submodule of $M$. 
\et   

Till end of this section, we assume that  Zhu's algebra $\Av$ is
semisimple.

\bl{L1} Let $M$ be an admissible $V$-module. For $\lm\in \C$, the
subspace $\ops_{n\in
\Z}M^{(\lm+n)}$ is an admissible $V$-submodule of $M$. Therefore, 
$N=\ops_{\lm\in\C}M^{(\lm)}$ is also an admissible $V$-submodule of M. 
\el
\proof Set
$W=\ops_{n\in\Z}M^{(\lm+n)}$. We have to prove that 
\eqa
\T{a}(n)W\subset W\hbox{ for all $a\in V$ and $n\in
\Z$}\label{AS1}
\eeqa
and that
\eqa d(W)\subset W,\label{AS2}
\eeqa
where $d$ is the degree operator of $M$. First we prove (\ref{AS1}).  For
$a\in V$ and $n\in
\Z$, we have a commutation relation 
\eqa [L(0),\T{a}(n)]=-n\,\T{a}(n).\label{CR1}
\eeqa
Hence for any $u\in M$,
\eqn (L(0)-\lm+n)\T{a}(n) u&=&L(0)\T{a}(n)u-(\lm-n)\T{a}(n)u\\
&=&\T{a}(n)L(0)u-n\T{a}(n)u-(\lm-n)\T{a}(n)u\\
&=&\T{a}(n)(L(0)-\lm)u.
\eeqn
This implies that if $u\in M^{(\lm)}$, then $\T{a}(n)u\in M^{(\lm-n)}$.
Thus we see that $\T{a}(n)u\in W$  for any $u\in W$.  Next we show
(\ref{AS2}). Let $u\in M^{(\lm)}$, and choose $k\in\Z_{+}$ so that
$(L(0)-\lm)^{k}u=0$. Since
$[d,L(0)]=0$ by (\ref{CR2}), we have 
$$(L(0)-\lm)^{k}du=d(L(0)-\lm)^{k}u=0.$$
Hence $du\in M^{(\lm)}$. This shows that $W$ is
invariant under the action of $d$. Thus we have (\ref{AS2}). \qed

Now we can prove the following proposition. 

\bp{MP1} Let $V$ be a vertex operator algebra, and suppose that Zhu's
algebra $A(V)$ is semisimple. Then any admissible $V$-module $M$
is decomposed into a direct sum of generalized eigenspaces for
$L(0)$, i.e., $M=\ops_{\lm\in\C}M^{(\lm)}.$ 
\ep
{\proof} Set $N=\ops_{\lm\in
\C}M^{(\lm)}$. By Lemma
\ref{L1}, $N$ is an admissible submodule of $M$. We need to prove
that $M/N=0$. Thus we show that $(M/N)_{n}=0$ for any
$n\in\N$  by induction on $n$. First we prove $(M/N)_{0}=0$. By
Proposition \ref{Pro1} (1),  $M_{0}$ is an
$\Av$-submodule of $\W(M)$. Since $\Av$ is semisimple, $L(0)$ acts on
$M_{0}$ diagonally. Hence $M_{0}$ is  contained in $N$. This implies that 
$(M/N)_{0}=0$. Let $n$ be a positive integer, and assume that
$(M/N)_{m}=0$ for any $m< n$. We prove $(M/N)_{n}=0$. By the
induction hypothesis, we see that  the subspace
$(M/N)_{n}$ is contained in $\W(M/N)$ and is invariant under the action
of $o(a)$ for any $a\in V$. Hence by Proposition \ref{Pro1} (2),
$(M/N)_{n}$ is an $\Av$-module, and 
$L(0)$ acts on $(M/N)_{n}$ diagonally. Therefore, it suffices to show that
 $u\in N$ for any $L(0)$-eigenvector $u+N\cap M_{n}\in (M/N)_{n}$. 

Let $u+N\cap M_{n}$ be an $L(0)$-eigenvector with eigenvalue $\lm$.
Then there exists $v\in N$ such that $L(0)u=\lm u +v$.
Since $v\in N$, we can write $v=\sum_{i=1}^{d}v_{i}$ for some $v_{i}\in
M^{(\lm_{i})}$. For $1\leq i\leq d$, choose $k_{i}\in\N$ so that 
$(L(0)-\lm_{i})^{k_{i}}v_{i}=0$, and set $f(x)=\prod_{i=1,\lm_{i}\neq
\lm}^{d}(x-\lm_{i})^{k_{i}}\in\C[x].$  Then it is obvious that 
\eqa
(L(0)-\lm)^{k}f(L(0))v=0\label{PL1}
\eeqa
 for some $k\in\N$. Thus it follows from (\ref{PL1}) that
\eqn (L(0)-\lm)^{k+1} f(L(0))u&=&f(L(0))(L(0)-\lm)^{k+1}u\\
&=&f(L(0))(L(0)-\lm)^{k}v\\
&=&0.
\eeqn
This implies that $f(L(0))u\in M^{(\lm)}\subset N$.  On the other
hand, we have $f(L(0))u-f(\lm)u\in N$  because $L(0)u-\lm u=v\in N$
and $N$ is $L(0)$-invariant. Since $f(\lm)$ is nonzero, we see that  $u\in
N$. This implies $(M/N)_{n}=0$.\qed

Let $M$ be an admissible $V$-module. Then by Lemma \ref{L1}, for
$\lm\in\C$, the subspace $\ops_{n\in\Z}M^{(\lm+n)}$ is an admissible
$V$-submodule of $M$. We next prove that there is an integer $m\in\Z$
such that $M^{(\lm+n)}=0$ for all integer $n<m$, more precisely we have: 

\bp{MP2} Let $M$ be an admissible $V$-module. If
$\lm\in\C-(\mathcal{P}(V)+\N)$, then $M^{(\lm)}=0$. 
\ep 
\proof Set $P'=\mathcal{P}(V)+\N$. For $\lm\in\C$, we have a
direct sum decomposition $M^{(\lm)}=\ops_{n=0}^{\infty}M^{(\lm)}\cap
M_{n}$, since $M^{(\lm)}$ is $d$-invariant. Set
$M_{n}^{(\lm)}=M^{(\lm)}\cap M_{n}$. By using
induction on $n$, we prove that $M_{n}^{(\lm)}=0$ for any
$\lm\in \C-P'$ and $n\in\N$. First we show that $M^{(\lm)}_{0}=0$ for
any $\lm\in\C-P'$. By Proposition \ref{Pro1} (1), $M_{0}^{(\lm)}$ is an
$\Av$-submodule of
$\W(M)$ for any $\lm\in\C$. The semisimplicity of $\Av$ implies that
$L(0)$  acts on
$M_{0}^{(\lm)}$ by the scalar $\lm$.  Thus by Proposition \ref{Pro1} (2),
we see that $M_{0}^{(\lm)}$ is zero if $\lm\in\C-P'$. Let
$n_{0}\in\Z_{+}$, and assume that
$M_{n}^{(\lm)}=0$ for all $\lm\in\C-P'$ and $n<n_{0}$. Then we show that
$M_{n_{0}}^{(\lm)}=0$ for any $\lm\in\C-P'$. By (\ref{CR2}) and
(\ref{CR1}), we have 
\eqa \T{a}(m)M_{n_{0}}^{(\lm)}\subset
M_{n_{0}-m}^{(\lm-m)}\label{Fact1}
\eeqa
for $a\in V$, $m\in\Z$ and $\lm\in\C$. If  $\lm\in\C-P'$
then $\lm-m\in\C-P'$ for any $m\geq1$. Hence by induction hypothesis,
we see that
$\T{a}(m)M_{n_{0}}^{(\lm)}=0$ for $m\geq1$ and $\lm\in\C-P'$. This
implies that $M_{n_{0}}^{(\lm)}\subset\W(M)$. By (\ref{Fact1}), the space
$M_{n_{0}}^{(\lm)}$ is invariant under the action of
$o(a)$ for all $a\in V$. Thus by Proposition
\ref{Pro1} (1), $M_{n_{0}}^{(\lm)}$ is an
$\Av$-submodule of $\W(M)$.  By the same reason as in the case $n=0$,
we have
$M_{n_{0}}^{(\lm)}=0$ for any $\lm\in\C-P'$. This completes the proof.
\qed     

Let us prove Theorem \ref{MT}. 

\noindent
{\it proof of Theorem \ref{MT}.} By Proposition \ref{MP1} and Proposition
\ref{MP2}, we  see that any admissible $V$-module $M$ is decomposed
as follows:
$$M=\ops_{\lm\in I}(\ops_{n\in\N}M^{(\lm+n)}),$$
where $I$ is an subset of $\Pv{V}$. Since $\Pv{V}$ is a
finite set by Theorem \ref{T2}, the subset $I$ is also finite. If we set
$I=\{\lm_{1},\ldots,\lm_{s}\}$, then we have the decomposition
$$M=\ops_{i=1}^{s}(\ops_{n\in\N}M^{(\lm_{i}+n)}).$$
By Lemma \ref{L1} and Proposition \ref{MP2}, we see that
$\ops_{n\in\N}M^{(\lm_{i}+n)}$ is admissible $V$-submodule of $M$ for
each $\lm_{i}$. This completes the proof.\qed

Let $V$ be a vertex operator algebra $V$ whose Zhu's algebra is
semisimple.  By Theorem \ref{MT}, we see that  if any admissible
 $V$-module of the form $ M=\ops_{n\in\N}M^{(\lm+n)}$ is completely
reducible for any $\lm\in \mathcal{P}(V)$, then $V$ is rational. Now we
prove the following proposition.

\bp{Pro2} If $\lm\in\Pv{V}$ satisfies a condition
$(\lm+\Z_{+})\cap\mathcal{P}(V)=\emptyset$, then any admissible
$V$-module of the form $M=\ops_{n\in\N}M^{(\lm+n)}$ is completely
reducible. 
\ep

To show this proposition, we first prove the following lemma. 


\bl{Lem3} Let $\lm$ and $M$ be as in Proposition \ref{Pro2}. If
$M^{(\lm)}=0$, then $M=0$. 
\el
\proof Assume that $M\neq0$. Then we can easily see that
$\W(M)$ is nonzero. By Proposition
\ref{Pro1} (2) and the assumption that
$(\lm+\Z_{+})\cap\Pv{V}=\emptyset$, we have
$\W(M)=M^{(\lm)}$. Therefore, $M^{(\lm)}$ is nonzero. \qed

As a corollary of Lemma \ref{Lem3}, we have:
\bc{Cor1}Let $\lm$ and $M$ be as in Proposition \ref{Pro2}. Then the
following hold{\rm :}

(1) Let $N$ be an admissible submodule of $M$. If $N$
contains $ M^{(\lm)}$, then $N=M$. 

(2) $M$ is irreducible as an admissible $V$-module if and
only if $M^{(\lm)}$ is irreducible as an $\Av$-module. 
\ec
\proof  It is clear that if $M^{(\lm)}\subset N$, then
$N^{(\lm)}=M^{(\lm)}$. Thus the quotient module
$M/N$ is zero by Lemma \ref{Lem3}. This proves (1). Next we prove
(2). As in the proof of Lemma \ref{Lem3}, we see that $\W(M)=M^{(\lm)}$.
Thus by Theorem
\ref{T2} (3) if
$M$ is an irreducible $V$-module, then
$M^{(\lm)}$ is irreducible as an $\Av$-module.
Conversely, assume that $M^{(\lm)}$ is irreducible as an $\Av$-module.
Let $N$ be a nonzero admissible submodule of $M$. By Lemma
\ref{Lem3},
$N^{(\lm)}$ is nonzero. Hence $N^{(\lm)}$ is a nonzero $\Av$-submodule
of $M^{(\lm)}$. Therefore, by the irreducibility of $M^{(\lm)}$,
$N^{(\lm)}=M^{(\lm)}$. Then (1) implies that $N=M$. \qed

By using Corollary \ref{Cor1}, we prove Proposition \ref{Pro2}. 

\noindent
{\it Proof of Proposition \ref{Pro2}.\ } Since $M^{(\lm)}$ is an
$\Av$-module, it is a direct sum of irreducible $\Av$-modules as
$$M^{(\lm)}=\ops_{i\in I}W_{i},$$
where $I$ is an index set  and $W_{i}$ is an irreducible component of
$M^{(\lm)}$. Let
$u_{i}\in W^{i}$ be a nonzero element, and set $N^{i}=\B u_{i}\K$. By
Lemma \ref{Lem1}, it is clear that 
$$(N^{i})^{(\lm+n)}=\span_{\C}\{\,\T{a}(-n)u_{i}\,|\, a\in V\,\}$$ for all
$n\in\N$. In particular, we see that $(N^{i})^{(\lm)}$ is a nonzero
$\Av$-submodule of $W^{i}$. Hence $(N^{i})^{(\lm)}=W^{i}$. Thus  it
follows from Corollary \ref{Cor1} (2) that $N^{i}$ is irreducible. Put
$N=\sum_{i\in I}N^{i}$. Then $N$ is completely reducible as an
admissible $V$-module. Since $N$ contains $M^{(\lm)}$,
Corollary \ref{Cor1} (1) shows that $N=M$. Thus $M$ is completely
reducible. \qed

By Theorem \ref{MT} and Proposition \ref{Pro2}, we have:
\bc{Thm3} Let $V$ be a vertex operator algebra. Suppose that Zhu's
algebra
$\Av$ is semisimple. If $(\lm+\Z_{+})\cap\mathcal{P}(V)=\emptyset$ for
any $\lm\in\Pv{V}$, then $V$ is rational. 
\ec

We shall give a remark on similar results for admissible
twisted modules. 

Let $V$ be a vertex operator algebra, and let $g$ be an
automorphism of $V$. In \cite{DLM1}, Dong, Li and Mason introduced an
associative algebra
$A_{g}(V)$,which is a generalization of  Zhu's
algebra $A(V)$. It is known that all statements in Theorem \ref{T2} hold
for $A_{g}(V)$ and
admissible $g$-twisted $V$-modules, instead of $A(V)$ and admissible
$V$-modules.  Since in the proofs of Theorem
\ref{MT} and Proposition
\ref{Pro2} we use only Theorem \ref{T2} and results which follow from
Theorem \ref{T2}, we see that the similar statements in Theorem
\ref{MT} and Proposition \ref{Pro2} hold for any $g$-twisted
$V$-modules:

\bt{MT2} Let $V$ be a vertex operator algebra, and let $g$ be an
automorphism of order $T$. Suppose that the associative algebra
$A_{g}(V)$ is semisimple. Then any admissible $g$-twisted $V$-module
$M$ is decomposed into a direct sum of generalized eigenspaces for $L(0)$
as follows:
$$M=\ops_{i=1}^{s}(\ops_{n\in\N}M^{(\lm_{i}+n/T)})$$
for some $\lm_{i}\in\mathcal{P}_{g}(V)$ $(i=1,2, \cdots,s)$,
where  $\mathcal{P}_{g}(V)$ is the set of the lowest weights of all
irreducible admissible $g$-twisted $V$-modules. Moreover, for each
$\lm_{i}$ the subspace $\ops_{n\in\N}M^{(\lm_{i}+n/T)}$ is an admissible
$g$-twisted submodule of $M$. 
\et
\bp{MP4} Let $V$, $g$ and $T$ be as in Theorem \ref{MT2}. If
$\lm\in\mathcal{P}_{g}(V)$ satisfies a condition
$(\lm+(1/T)\Z_{+})\cap\mathcal{P}_{g}(V)=\emptyset$,  then any
admissible $g$-twisted $V$-module of the form
$M=\ops_{n\in\N}M^{(\lm+n/T)}$ is completely reducible. 
\ep 

Therefore, as a corollary of Proposition \ref{MP4} we have:
\bc{Cor4} Let $V$, $g$ and $T$ be as in Theorem \ref{MT2}. If the
condition 
$(\lm+(1/T)\Z_{+})\cap\mathcal{P}_{g}(V)=\emptyset$ holds for any
$\lm\in\mathcal{P}_{g}(V)$, then $V$ is $g$-rational. 
\ec

\section{Main Results}\label{S4.0}
In this section we prove that the vertex operator algebra
$\V{+}$ for a rank one even lattice $L=\Z\al$ is rational 
if $k=\B\al,\al\K/2$ is a prime integer. In Section \ref{S4.1}, we recall the
construction of the vertex operator algebra $\V{+}$ and its irreducible
modules. In Section
\ref{S4.2}, we show the rationality of
$\V{+}$ for a prime integer $k$ by using the results in the previous
section. 

\subsection{The Vertex Operator Algebra $\V{+}$}\label{S4.1}
Let $L=\Z\al$ be a rank one even lattice with a ${\Z}$-bilinear
form $\B\cdotp,\cdotp\K$ defined by $\B\al,\al\K=2k$ for a positive
integer $k$. Set 
$\h=\C\ots_{\Z} L$ and extend the form
$\B\cdotp,\cdotp\K$ to a $\C$-bilinear form on $\h$. 
Let $L^{\circ}$ be the dual lattice of $L$, and let
$\C[L^{\circ}]=\ops_{\lm\in L^{\circ}}\C e_{\lm}$ be the group algebra of
$L^{\circ}$. For a subset $M$ of $L^{\circ}$, we set
$\C[M]=\ops_{\lm\in M}\C e_{\lm}$. 

Let $\hh =\h\ots \C[t,t^{-1}]\ops\C K$ be a Lie algebra with
the commutation relations given by 
$$[X\ots t^{m},X'\ots
t^{n}]=m\,\del_{m+n,0}\B X,X'\K K\hbox{ and }[K,\hh]=0,\
(X,X'\in\h,m,n\in\Z).$$ 
Then $\hh^{+}=\h\ots
\C[t]\ops \C K$ is a subalgebra of $\hh$, and the group algebra
$\C[L^{\circ}]$ becomes a $\hh^{+}$-module by the actions 
$$(X\ots t^{n}).
e_{\lm}=\del_{n,0}\B X, \lm\K e_{\lm}\hbox{ and
}K.e_{\lm}=e_{\lm}\ (\lm\in L^{\circ},X\in \h,n\in\N).$$
 For a subset
$M$ of $L^{\circ}$  and $\mu\in L^{\circ}$, we set 
$$ V_{M}=U(\hh)\ots_{U(\hh^{+})}
\C[M]\hbox{ and }\Ml{\beta}=U(\hh)\ots_{U(\hh^{+})} \C e_{\beta},$$
where $U(\g)$ denotes the universal enveloping algebra of a Lie algebra
$\g$. Let $\pi$ be the representation of $\hh$ on
$V_{L^{\circ}}$, and $X(n)=\pi(X\ots t^{n})$ ($X\in\h,
n\in\Z$). We set $\1=1\ots e_{0}$ and $\w=(1/2)h(-1)^{2}e_{0}$, where
$h$ is a vector of $\h$ with $\B h,h\K=1$. Then it is known that there
exists a linear map
$Y:\V{}\to V_{L^{\circ}}[[z,z^{-1}]]$ such that
$(\V{},Y,\1,\w)$ has a simple vertex operator algebra structure, and
$(\Vl{\lm},Y)$ becomes an irreducible $\V{}$-module for $\lm\in
L^{\circ}$ (see \cite{FLM}). We remark that $(\Ml{0},Y,\1,\w)$ is an vertex
operator subalgebra of $\V{}$. We write $\M{}=\Ml{0}$ for simplicity. The
vertex operator $Y(e_{\beta},z)$ for $\beta\in L$ on
$V_{L^{\circ}}$ is given by
\eqa
Y(e_{\beta},z)=\exp\left(\sum_{n=1}^{\infty}
{\frac{\beta(-n)}{n}}z^{n}\right)
\exp\left(-\sum_{n=1}^{\infty}{\frac{\beta(n)}{n}}z^{-n}\right)
e_{\beta}z^{\beta(0)},\label{O1}
\eeqa 
where $e_{\beta}$ in the right-hand side means the left
multiplication of $e_{\beta}$ on the group algebra
$\C[L^{\circ}]$, and $z^{\beta(0)}$ is an operator on $V_{L^{\circ}}$ defined
by $z^{\lm(0)}u=z^{\B\lm,\mu\K}u$ for $\mu\in L^{\circ}$ and
$u\in\Ml{\mu}$. The operator $L(0)=\w_{1}$ is given by
$$L(0)=\frac{1}{2}h(0)^{2}+\sum_{n\in\Z_{+}}h(-n)h(n).$$
It is known that
$\V{}$ is a rational vertex operator algebra and every irreducible
$\V{}$-module is isomorphic to $\Vl{\lm}$ for some $\lm\in L^{\circ}$
(see \cite{D2}).  

Let $\th$ be a linear isomorphism of $V_{L^{\circ}}$ defined by 
\eqn 
\th(X_{1}(-n_{1})X_{2}(-n_{2})\cdots X_{\ell}(-n_{\ell})\ots
e_{\lm})=(-1)^{\ell}\,X_{1}(-n_{1})X_{2}(-n_{2})\cdots 
X_{\ell}(-n_{\ell})\ots e_{-\lm}
\eeqn
for $X_{i}\in\h,n\in\Z_{+}$ and $\lm\in L^{\circ}$. Then $\th$
induces an automorphism of $\V{}$ of order $2$. For a $\th$-invariant
subspace $W$ of $V_{L^{\circ}}$, we denote by $W^{\pm}$ the $(\pm
1)$-eigenspaces of $W$ for $\th$. Then $(\V{+},Y,\1,\w)$ is a simple
vertex operator algebra, and $\V{\pm}$, $V_{{\al}/{2}+L}^{\pm}$ and
$V_{r\al/2k+L}$ for
$1\leq r\leq k-1$ become irreducible $\V{+}$-modules (see \cite{DN2}). 

Another irreducible $\V{+}$-modules come from $\th$-twisted
$\V{}$-modules.  We briefly review the
constructions of $\th$-twisted
$\V{}$-modules (see \cite{FLM} and \cite{D2} for more detailed 
construction). Let
$\hh[-1]=\h\ots t^{1/2}\C[t,t^{-1}]\ops \C K$ be a Lie algebra with the 
commutation relation given by 
$$[X\ots t^{m},X'\ots
t^{n}]=m\,\del_{m+n,0}\,\B X,X'\K\, K\hbox{ and }[K,\hh[-1]]=0\
(X,X'\in\h,m,n\in \frac{1}{2}+\Z).$$
Then $\hh[-1]^{+}=\h\ots t^{1/2}\C[t]\ops \C K$ is a subalgebra of
$\hh[-1]$. Let $T_{1}$ and $T_{2}$ be irreducible $\C[L]$-modules
on which $e_{\al}$ acts as $1$ and $-1$ respectively. We define the action
of $\hh[-1]^{+}$ on $T_{i}$ by the actions
 $$(X\ots t^{1/2+n}).u=0\hbox{  and }K.u=u\ (X\in\h,n\in\N,u\in T_{i}).$$
Set 
$$ \V{T_{i}}=U(\hh[-1])\ots_{U(\hh[-1]^{+})}T_{i}\ (i=1,2),$$  then
$\V{T_{i}}$ has the irreducible $\th$-twisted $\V{}$-module structure. 

Now we define the action of the automorphism $\th$ on $\V{T_{i}}$
$(i=1,2)$ by
\eqn \th(X_{1}(-n_{1})\cdots
X_{\ell}(-n_{\ell})u)=(-1)^{\ell}X_{1}(-n_{1})\cdots X_{\ell}(-n_{\ell})u,
\eeqn 
for $X_{i}\in\h, n_{i}\in 1/2+\N$ and $u\in T_{i}$, and set $\V{T_{i},\pm}$
the $\pm 1$-eigenspaces of $\V{T_{i}}$ for $\th$. Then $\V{T_{i},\pm}$
($i=1,2$) become irreducible $\V{+}$-modules (see
\cite{DN2}). 

All irreducible $\V{+}$-modules are
completely classified in \cite{DN2}. Their result is as follows:


\bt{T1} The set
\eqa\{\V{\pm},\Va{\pm},\V{T_{i},\pm},\Vl{r\al/2k}\,
|\,i=1,2,1\leq r\leq k-1\,\}\label{IM1}
\eeqa
is the complete list of all inequivalent irreducible admissible
$\V{+}$-modules.
\et

The following result is also found in \cite{DN2};
\bp{Pro11} Zhu's algebra $A(\V{+})$ is semisimple.
\ep

It is known that $\V{+}$ is rational when $k=1,2$. In fact, if $k=1$, then
$\V{+}$ is isomorphic to the lattice vertex operator algebra $V_{L'}$,
where $L'=\Z\beta$ is an even lattice with
$\B\beta,\beta\K=8$ (see \cite{DG}). If $k=2$, then $\V{+}$ is isomorphic
to the vertex operator algebra $L(1/2,0)\ots L(1/2,0)$, where
$L(1/2,0)$ is the minimal vertex operator algebra of central charge
$1/2$ (see \cite{DGH}). 

\subsection{Proof of the rationality of $\V{+}$ for a prime integer $k$
}\label{S4.2}
First we consider the set of lowest weights of $\V{+}$. By the
construction of irreducible $\V{+}$-modules, we have the following table. 
\begin{center}
\begin{tabular}{|c|c|c|c|c|c|}\hline
$\V{+}$&$\V{-}$&$\Vl{\lm_{r}}$ $(1\leq r\leq
k-1)$&$\Va{\pm}$&$\V{T_{1},+},\V{T_{2},+}$&$\V{T_{1},-},\V{T_{2},-}$\\  \hline
$0$&$1$&$r^{2}/4k$&$k/4$&$1/16$&$9/16$\\  \hline
\end{tabular}
\end{center}
\begin{center}
{{\bf Table 1.} The lowest weights of irreducible $\V{+}$-modules.}
\end{center}
 Hence we have
\eqa \mathcal{P}(\V{+})=\{\,0,1,r^{2}/4k,1/16,9/16\,|\,1\leq
r\leq k\,\}.\label{Wt1}
\eeqa
Now we shall apply the results in Section \ref{S3} to the case that
$V=\V{+}$. Since $0,1\in\Pv{\V{+}}$, we can
not use Theorem \ref{Thm3} to prove the rationality of $\V{+}$. But in the
case
$k=\B\al,\al\K/2$ is a prime integer, we have the following lemma:


\bl{Lem4} Suppose that $k$ is a prime integer. 
Then the lowest weights $0,1,1/16,9/16$ and $r^{2}/4k$
for $1\leq r\leq k$ are all distinct. Furthermore, 
for any $\lm\in\Pv{\V{+}}-\{0\}$,
we have $(\lm+\Z_{+})\cap\Pv{\V{+}}=\emptyset$.
\el
\proof We have to prove
that $r^{2}/4k, r^{2}/4k-1/16,r^{2}/4k-9/16\not\in\Z$ for $1\leq r\leq k$
and that 
$(r^{2}-s^{2})/4k\not\in\Z$ for any $1\leq r\neq s\leq k$.
We first show that
$r^{2}/4k$ is not integer for $1\leq r\leq k$. If $r^{2}/4k$ is an integer,
then $k|r$ because
$k$ is prime. So $r=k$ and then $r^{2}/4k=k/4$ is an integer. This is
contradiction. Next we assume that $r^{2}/4k-1/16$ is an integer. Then
$4r^{2}-k\equiv 0\pmod{16}$. Hence $k\equiv 0 \pmod{4}$. This is
contradiction. We can also prove that
$r^{2}/4k-9/16$ is not integer in the same way. Finally we show that
$(r^{2}-s^{2})/4k$ is not in $\Z$ if $r\neq s$. Assume that
$(r^{2}-s^{2})/4k$ is an integer for some $1\leq r< s\leq k$. Then $k|(r+s)$
since
$k$ is prime and $1\leq s-r\leq k-1$. Furthermore since $3\leq r+s\leq
2k-1$, we have  $r+s=k$ and $k\geq 3$.   Hence $(r^{2}-s^{2})/4k
=(2s-k)/4$.  It contradicts with the fact that $k$ is an odd prime integer.
\qed

Thus we see that if $k$ is prime then any nonzero lowest weights for
$\V{+}$ satisfies the condition in Proposition \ref{Pro2}. Therefore, we
have: 


\bp{Pro3} If $k$ is prime and $\lm\in\Pv{\V{+}}$ is nonzero, then any
admissible $\V{+}$-module of the form
$M=\ops_{n=0}^{\infty}M^{(\lm+n)}$ is completely reducible. In particular,
if $M^{(\lm)}$ is irreducible as an $A(\V{+})$-module, then $M$ is
irreducible as an admissible $\V{+}$-module.
\ep
{\bf Remark.}  If $k$ is not prime, then there exist integers $r,s$ with
$0\leq r< s\leq k$ such that $r^{2}/4k-s^{2}/4k\in\Z$. In fact, when 
$k$ is not prime, $k$ can be written as $k=pqn$ for two positive prime
integers
$p,q$ ($q\leq p$) and $n\in\Z_{+}$. Then we see that  $0\leq
np-nq < np+nq \leq npq=k$, and that
$(np+nq)^{2}/4k-(np-nq)^{2}/4k=n\in\Z$.

{\bf Here and further we suppose that $k$ is prime}.  
By Proposition \ref{Pro3} and Theorem \ref{MT}, to show that $\V{+}$ is
rational, it suffices to prove that every admissible
$\V{+}$-module $M$ of the form $M=\ops_{n=0}^{\infty}M^{(n)}$ is  
completely reducible.

Let $W$ be the admissible
$V$-submodule generated by $M^{(0)}$ and $U$ the admissible
$V$-submodule generated by $M^{(1)}$. Since the sum $W+U$ contains 
$M^{(0)}\ops M^{(1)}$, the quotient module $M/(W+U)$ is decomposed as
$M/(W+U)=\ops_{n=2}^{\infty}(M/(W+U))^{(n)}$. Thus by Lemma
\ref{Lem3} and Lemma \ref{Lem4}, we see that $M/(W+U)=0$, that is,
$M=W+U$. Therefore, to prove that $M$ is completely reducible, it is
enough to show that $U$ and $W$ are completely reducible.
 
We first show that the admissible $\V{+}$-submodule $W$ is a direct sum
of copies of $\V{+}$. To prove this we need the following proposition:


\bp{Pro5} For any admissible $\V{+}$-module of the form 
$M=\ops_{n=0}^{\infty}M^{(n)}$, we have $L(-1)M^{(0)}=0$. 
\ep 
\proof Let $u\in M^{(0)}$, and suppose that $L(-1)u\neq 0$. 
We see that for any homogeneous vector $a\in\V{+}$ with positive weight, $o(a)u=0$ since $M^{(0)}$ is a direct sum of copies of the irreducible
$A(\V{+})$-module $\W(\V{+})=\C\1$. Hence we have
$\widetilde{(L(-1)a)}(n)u=0$ for any $a\in\V{+}$ and $n\geq0$. This
shows that for $a\in\V{+}$ and
$n\geq0$, 
\eqn \T{a}(n+1)L(-1)u&=&L(-1)\T{a}(n+1)u-\widetilde{(L(-1)a)}(n)u\\
&=&0.
\eeqn
Thus $L(-1)u\in\W(M)\cap M^{(1)}$. We note that the subspace
$\W(M)\cap M^{(1)}$ is an $A(\V{+})$-module on which $L(0)$ acts by the scalar $1$. 
By Table 1 and Lemma \ref{Lem4}, the subspace $\W(M)\cap
M^{(1)}$ is a direct sum of
copies of the irreducible $A(\V{+})$-module $\W(\V{-})=\C\al(-1)\1$.
This implies that $\C L(-1)u\cong \C\al(-1)\1$ as $A(\V{+})$-modules.
Set $N=\B L(-1)u\K$. Then we see that $N$ is an
admissible $V$-submodule of $M$ such that
$N^{(0)}=0$ and $N^{(1)}=\C L(-1)u$.  Hence by Proposition \ref{Pro3},
$N$ is an irreducible $\V{+}$-module isomorphic to $\V{-}$. 
Put 
\eqa E=e_{\al}+e_{-\al}\in \V{+}\hbox{ and }F=e_{\al}-e_{-\al}\in
\V{-}.\label{EF1}
\eeqa
 By
direct calculations, we have 
$E(0)\al(-1)\1=-2kF$ and $E(n)\al(-1)\1=0$ for all $n\geq1$. Hence we
see that
\eqa
E(0)L(-1)u\neq0\hbox{ and }E(n)L(-1)u=0\hbox{ for all
$n\geq1$}.\label{FG1}
\eeqa
 Since $[L(-1),E(n)]u=-nE(n-1)u$ for any $n\in\Z$,
$L(-1)E(n)u=-nE(n-1)u$ for $n\geq1$. Thus $E(k-1)u=0$ implies that
$E(0)u=0$. Then
$E(0)L(-1)u=L(-1)E(0)u=0$. This contradicts with (\ref{FG1}).
Therefore we conclude that $L(-1)u=0$.  \qed

Let $u$ be a nonzero vector in $M^{(0)}$. Then $\B u\K$ is an admissible
$\V{+}$-submodule  of $M$ such that $\B u\K^{(0)}=\C u$.  Now we show
that $\B u\K^{(1)}=0$. Since $\B u\K^{(1)}$ is spanned by vectors of the
form $a(\wt(a)-2)u$ for homogeneous vectors $a\in
\V{+}$, it is enough to prove that $a(\wt(a)-2)u=0$ for
any homogeneous vector $a\in \V{+}$. We may assume that
$\wt(a)\neq1$ because $(\V{+})_{1}=0$ unless $k=1$. Then by Proposition
\ref{Pro5} we have $a(\wt(a)-2)u=-(\wt(a)-1)^{-1}[L(-1),o(a)]u=0$. 
Thus $\B u\K^{(1)}=0$. Therefore we have 
$\B u\K=\C u\ops (\ops_{n=2}^{\infty}\B u\K^{(n)})$. Let $N$ be a
nonzero admissible $\V{+}$-submodule of $\B u\K$. By Proposition
\ref{Pro1} (2) and Lemma \ref{Lem4}, we see that $N^{(0)}\neq0$. This
implies that $N=\B u\K$. Hence $\B u\K$ is an irreducible
$\V{+}$-module, and is isomorphic to $\V{+}$ by Table 1. Hence we have
the following proposition:


\bp{Pro6} The admissible $\V{+}$-submodule $W$  is a sum of
irreducible submodules isomorphic to $\V{+}$. Thus $W$ is completely
reducible. In particular, $W^{(1)}$ is zero. 
\ep

Next we show that the admissible $V$-submodule $U$, which is the
admissible $\V{+}$-submodule of $M$ generated by $M^{(1)}$,  is a sum of
copies of the irreducible
$\V{+}$-module
$\V{-}$. To prove this we have to show the following proposition.


\bp{Pro7} For any admissible $\V{+}$-module of the form 
$M=\ops_{n=0}^{\infty}M^{(n)}$, we have $L(1)M^{(1)}=0$. 
\ep

Before giving the proof of Proposition \ref{Pro7}, we prove some lemmas.


\bl{Lem5}(1) If $a\in \Ml{0}^{+}$ is homogeneous, then 
$$Y(E,z)a\in z^{-\wt(a)}(\V{+})[[z]].$$
(2) For any homogeneous vector $a\in\V{+}$ with $\wt(a)<k$, we have
$E(2k-2)a=0$. 
\el
\proof First we show the assertion (1). Let $a$ be a homogeneous vector
in $\M{}$. Then we have 
\eqa Y(e_{\al},z)a=\exp\left(\sum_{n=1}^{\infty}
{\frac{\al(-n)}{n}}z^{n}\right)
\exp\left(-\sum_{n=1}^{\infty}{\frac{\al(n)}{n}}z^{-n}\right)
a\ots e_{\al}.\label{Act1}
\eeqa
Now set  
$$\exp\left(-\sum_{n=1}^{\infty}{\frac{\al(n)}{n}}z^{-n}\right)
=\sum_{n=0}^{\infty}p_{n}(\al)z^{-n}.$$
Then $p_{n}(\al) a\ots e_{\al}\in \Ml{\al}$ for $n\geq0$. Since the weight
of $p_{n}(\al) a\ots e_{\al}$ is $\wt(a)-n+\wt(e_{\al})$, we see that  if
$n>\wt(a)$, then $p_{n}(\al) a\ots e_{\al}=0$. Thus we have 
$$\exp\left(-\sum_{n=1}^{\infty}{\frac{\al(n)}{n}}z^{-n}\right)a\ots
e_{\al}\in z^{-\wt(a)}(\V{})[z],$$
hence $Y(e_{\al},z)a\in z^{-\wt(a)}(\V{})[[z]]$. By the same way, we have 
$Y(e_{-\al},z)a\in z^{-\wt(a)}(\V{})[[z]]$. Therefore if $a\in\M{+}$ is
homogeneous, then $Y(E,z)a\in z^{-\wt(a)}(\V{+})[[z]]$. This proves (1).
Next we show the assertion (2). It is clear that for a homogeneous vector
$a\in\V{+}$, if $\wt(a)< k$, then $a\in\M{+}$. Hence by (1), $E(n)a=0$ for
$n\geq\wt (a)$. In particular, we see that
$E(k-1)a=0$.\qed 


\bl{Lem6}  For any $v\in W^{(k-1)}$, $E(2k-2)v=0$. 
\el
\proof By Proposition \ref{Pro6}, $W$ is a direct sum of copies of the
irreducible $\V{+}$-module $\V{+}$. Thus $W^{(k-1)}$ is a direct sum of
$(\V{+})_{k-1}$. Hence by Lemma
\ref{Lem5} (2), $E(2k-2)u=0$ for any
$u\in W^{(k-1)}$. \qed

{\it Proof of Proposition \ref{Pro7}.} Consider the
quotient module $M/W$. Since $M^{(0)}=W^{(0)}$, we have 
$M/W=\ops_{n=1}^{\infty}(M/W)^{(n)}$. Hence by Proposition \ref{Pro3},
$M/W$ is completely reducible and is a direct sum of copies of the
irreducible $\V{+}$-module $\V{-}$. Since $W^{(1)}$ vanishes, $M^{(1)}$
is embedded in $M/W$. So a nonzero $u\in M^{(1)}$ generates an
irreducible $\V{+}$-submodule of $M/W$ isomorphic to $\V{-}$ such
that
$\C u\cong \C\al(-1)\1$ as $A(\V{+})$-modules. Since $E(0)\al(-1)\1=-2k
F$ and $E(2k-2)F=-2\al(-1)\1$, we have
$E(2k-2)E(0)\al(-1)\1=4k\al(-1)\1$. We have also
$E(1)\al(-1)\1=0$. These imply that 
$E(2k-2)E(0)u=4ku$ and $E(1)u\in W$. Therefore we have 
\eqa
4kL(1)u&=&L(1)E(2k-2)E(0)u\nn\\
&=&E(2k-2)L(1)E(0)u\nn\\
&=&E(2k-2)(E(0)L(1)u+2(k-1)E(1)u).\label{FG2}
\eeqa
Then the fact that $E(0)L(1)u, E(1)u\in W^{(k-1)}$ and Lemma \ref{Lem6}
prove that the right hand side in (\ref{FG2}) is zero. Hence $L(1)u=0$.
\qed

Now suppose that $u\in M^{(1)}$. Let $a\in \V{+}$ be a quasi-primary
vector, i.e., a homogeneous vector which satisfies $L(1)a=0$. It is clear that
$\T{a}(n)u=0$ for $n\geq2$. By the commutation
relation (\ref{CR}) and (\ref{DP1}) we have $[L(1),
o(a)]u=(\wt(a)-1)\T{a}(1)u$. Thus Proposition \ref{Pro6} shows that
$(\wt(a)-1)\T{a}(1)u=0$. Since $(\V{+})_{1}=0$, we see that 
$\T{a}(1)u=0$. Hence $\T{a}(n)u=0$ for any quasi-primary vector
$a\in\V{+}$ and $n\in\Z_{+}$. By
$L(-1)$-derivative property (\ref{DP1}), we have 
$$\widetilde{(L(-1)^{m}a)}(n+1)u
=(-1)^{m}m!\Com{\wt(a)+m+n}{m}\T{a}(n+1)u=0$$
for $m,n\in\N$. On the other hand, since $(\V{+})_{1}=0$, we see that 
$\V{+}$ is spanned by vectors of the form $L(-1)^{m}a$ for $m\in\N$ and
quasi-primary vectors $a\in\V{+}$. Thus we see that $\T{a}(n)u=0$
for any $a\in\V{+}$ and
$n\in\Z_{+}$. This implies that $u\in\W(M)$. Therefore, we have
$M^{(1)}\subset \W(M)$. By Lemma \ref{Lem1}, we have $U^{(0)}=0$. 
Therefore it follows from Proposition \ref{Pro3} that  $U$ is completely
reducible. Furthermore by Lemma \ref{Lem4} and Table 1, we have:

\bp{Pro9} The admissible $V$-submodule $U$ is a direct sum of copies of
the irreducible
$\V{+}$-module $\V{-}$. 
\ep 

Consequently, by Proposition \ref{Pro6} and Proposition \ref{Pro9}, we
see that $M=W+U$ is completely reducible. Thus we have:


\bp{Pro10} Any admissible $\V{+}$-module
$M=\ops_{n=0}^{\infty}M^{(n)}$ is completely reducible.
\ep

The following theorem follows from Proposition \ref{Pro3}, Proposition
\ref{Pro10} and Theorem \ref{MT}: 
\bt{Th5} If $k=\B\al,\al\K/2$ is a prime integer, then the vertex
operator algebra $\V{+}$ is rational.
\et

{\bf Acknowledgments:} I would like to thank Professor Kiyokazu
Nagatomo and Dr. Yoshiyuki Koga for useful suggestions, and I also
thank Professor Ching Hung Lam for reading this manuscript and giving
significant comments.


\begin{thebibliography}{100000}

\bibitem[A]{A}
T. Abe, Fusion rules for the charge conjugation orbifold, 
math.QA./0006101.

\bibitem[D1]{D1}
C. Dong, Vertex algebras associated with even lattices, {\it J. Algebra}
{\bf160} (1993), 245-265.

\bibitem[D2]{D2}
C. Dong, Twisted modules for vertex algebras associated with even
lattices, {\it J. Algebra} {\bf165} (1994), 90-112.

\bibitem[DG]{DG}
C. Dong and R. L. Griess, Rank one lattice type vertex operator
algebras and their automorphism groups, {\it J. Algebra} {\bf208}
(1998), 262-275.

\bibitem[DGH]{DGH}
C. Dong, R. L. Griess and G. H{\"o}hn, Framed vertex operator algebras,
codes and the moonshine module, {\it Commun. Math. Phys.} {\bf193}
(1998), 407-448.

\bibitem[DL]{DL}
C. Dong and J. Lepowsky, ``Generalized vertex algebras 
and relative vertex operators'', Progress in Math., Vol.112, Birkh{\"a}user,
Boston, 1993.

\bibitem[DLi]{DLi}
C. Dong and  Z. Lin, Induced modules for vertex operator
algebras, {\it Commun. Math. Phys.} {\bf179} (1996), 157-184.

\bibitem[DLM1]{DLM2}
C. Dong, H.-S. Li and G. Mason, Regularity of rational vertex
operator algebras, {\it Adv. in Math.} {\bf132} (1997), 148-166.

\bibitem[DLM2]{DLM1}
C. Dong, H.-S. Li and G. Mason, Twisted representations of vertex
operator algebras, {\it Math. Ann.} {\bf310} (1998), 571-600.

\bibitem[DM]{DM}
C. Dong and G. Mason, On quantum Galois theory, {\it Duke Math. J.} {\bf86}, No. 2
(1997), 305-321.

\bibitem[DMZ]{DMZ}
C. Dong, G. Mason and Y.-C. Zhu,  Discrete series of the Virasoro
algebra and the moonshine module, in ``Proc. Sympos. Pure Math. Vol.
56," Chap. II,  pp. 295-316, {\it Amer. Math. Soc.}, Providence, 1994.

\bibitem[DN]{DN2}
C. Dong and K. Nagatomo, Representations of Vertex operator algebra
$\V{+}$ for rank one lattice $L$, {\it Commun. Math. Phys.} {\bf202}
(1999), 169-195.

\bibitem[DVVV]{DVVV}
R. Dijkgraaf, C. Vafa, E. Verlinde and H. Verlinde, The operator
algebra of orbifold models, {\it Commun. Math. Phys.} {\bf123} (1989),
485-526.

\bibitem[FHL]{FHL} 
I. Frenkel, Y.-Z. Huang and J. Lepowsky, On axiomatic approaches to
vertex operator algebras and modules, Mem. Amer. Math.
Soc. {\bf104} (1993).

\bibitem[FLM]{FLM}
I. Frenkel, J. Lepowsky and A. Meurman, ``Vertex Operator Algebras and
the Monster", Pure and Appl. Math., Vol. 134, Academic Press, Boston,
1988.

\bibitem[FZ]{FZ}
I. Frenkel and Y.-C. Zhu, Vertex operator algebras associated to
representations of affine and Virasoro algebras, {\it Duke Math. J.} {\bf66}
(1992), 123-168.

\bibitem[KT]{KT}
V. Kac and I. Todorov, Affine orbifolds and rational conformal field theory
extensions of $W_{1+\infty}$, 
{\it Commun. Math. Phys.} {\bf190} (1997), 57-111.

\bibitem[LS]{LS}
C.-H. Lam and H. Shimakura, {\it private communications}.

\bibitem[Li]{L1}
H.-S. Li, The regular representation, Zhu's $A(V)$-theory and induced
modules, math.QA/9909007.

\bibitem[M1]{M1}
M. Miyamoto, Griess algebras and conformal vectors in vertex operator
algebras, {\it J. Algebra} {\bf179} (1996), 523-548.

\bibitem[M2]{M2}
M. Miyamoto, Binary codes and vertex operator (super)algebras, {\it J.
Algebra} {\bf181} (1996), 207-222.

\bibitem[M3]{M3}
M. Miyamoto, Representation theory of code vertex operator algebras,
{\it J. Algebra} {\bf201} (1998), 115-150.

\bibitem[M4]{M5}
M. Miyamoto, A new construction of the Moonshine vertex operator
algebras over the real number field, {q-alg/9701012}.

\bibitem[Sh]{Sh}
H. Shimakura, Decomposition of the moonshine module with respect to a
code over $\Z_{6}$, master thesis, The University of Tokyo.

\bibitem[Z]{Z}
Y.-C. Zhu, Modular invariance of characters of vertex operator
algebras, {\it J. Amer. Math. Soc.} {\bf9} (1996), 237-302.

\end{thebibliography}
\end{document}